\documentclass[12pt]{article}
\usepackage{amsfonts,amsmath,amssymb,amscd,amsthm}
 \usepackage[russian]{babel}
\def\Z{{\mathbb Z}} \def\R{{\mathbb R}}

\long\def\comment#1\endcomment{}

\def\Z{{\mathbb Z}} \def\R{{\mathbb R}}

\long\def\comment#1\endcomment{}

\newtheoremstyle{mydefinition}
  {3pt}
  {3pt}
  {\normalfont}
  {\parindent}
  {\bfseries}
  {.}
  { }
  {}

\theoremstyle{mydefinition}
\newtheorem{pr}{}

\begin{document}


\comment

Existence proofs in combinatorics using independence

D. Ilyinskiy, A. Raigorodskiy and A. Skopenkov

This note is purely expository and is in Russian.
We show how to prove interesting combinatorial results using the {\it local Lovasz lemma}.
The note is accessible for students having basic knowledge of combinatorics;
the notion of independence is defined and the Lovasz lemma is stated and proved.
Our exposition follows `Probabilistic methods' of N. Alon and J. Spencer.
The main difference is that we show how the proof could have been invented.
The material is presented as a sequence of problems, which is peculiar not only to Zen monasteries but
also to advanced mathematical education; most problems are presented with hints or solutions.

\endcomment

\centerline{\uppercase{\bf Независимость и доказательства существования}}
\bigskip
\centerline{\uppercase{\bf в комбинаторике}
\footnote{
Заметка основана на занятиях, проведенных авторами на ФИВТ МФТИ,
в Кировской ЛМШ, Московской ВШ, а также на кружке `Олимпиады и математика'.
Благодарим А. Дайняка, К. Матвеева, А. Шаповалова, членов редколлегии сборника `Математическое Просвещение'
и участников занятий за полезные замечания и обсуждения.}
}

\bigskip
\centerline{\bf Д. Ильинский, А. Райгородский и А. Скопенков
\footnote{Все авторы: Московский Физико-Технический Институт.
А. Райгородский: Московский Государственный Университет.
А. Скопенков: Независимый Московский Университет.
Личные страницы: www.mccme.ru/\~{ }skopenko,
http://dm.fizteh.ru/staff.
А. Скопенков частично поддержан грантом фонда Саймонса.}
}

\bigskip
{\bf Введение. }

Цель этой заметки --- продемонстрировать метод доказательства некоторых интересных комбинаторных результатов  (пункты (b) задач \ref{prob-lovac-2cveta}-\ref{prolov-rainbow} и задачи \ref{prolov-hyp}-\ref{prolov-gen}), заключающийся в применении локальной леммы Ловаса \ref{prolov-sym}.
Для изучения заметки не нужно предварительных знаний, все необходимые понятия вводятся по ходу изложения.

Следующие две части введения важны, но формально не используются далее.

\smallskip
{\it Об открытии леммы Ловаса и ее роли в математике.}
Локальная лемма Ловаса была доказана в 1973 году выдающимся венгерским математиком Ласло Ловасом.
Впрочем, тогда Ловасу было всего 25 лет, и, хотя яркие результаты у него уже к тому времени были, все-таки на тот момент его воспринимали не как классика, но как восходящую звезду.
Он уже был трехкратным победителем международных математических олимпиад (1964, 1965 и 1966 годов).
Классиком Ловас станет позже, и весьма серьезную роль в этом сыграет доказанная им Локальная лемма.
Разумеется, не только она: будет и топологический метод в комбинаторике, и мощные результаты в теории алгоритмов, и значительный вклад в науку о графовых пределах, и многое другое.
Тем не менее, Локальная лемма --- это замечательный инструмент вероятностной комбинаторики, благодаря которому были получены и продолжают получаться многочисленные яркие результаты в области дискретной математики и теории алгоритмов.

Работа, в которой Ловас формулирует и доказывает свою Локальную лемму, написана в соавторстве с Полом Эрдешем --- еще одним великим специалистом по комбинаторике, основателем большой научной школы, автором множества задач и идей.
Среди прочего, Эрдеш был одним из самых активных пропагандистов вероятностного метода в комбинаторике.
Поэтому, несмотря на то, что Локальную лемму доказал именно Ловас, роль Эрдеша во всем этом
не стоит недооценивать.
В статье Эрдеша и Ловаса \cite{EL} речь шла о раскрасках {\it гиперграфов} (т.е. наборов подмножеств конечного множества).
Как раз ради доказательства существования некоторой раскраски Локальная лемма и придумывалась (т.е. ради обобщения задач \ref{prob-lovac-2cveta},  \ref{prolov-hyp} и \ref{prolov-subsets}; не бойтесь, они формулируются и решаются без слова `гиперграф').
Однако очень быстро стало ясно, насколько это мощный и плодотворный инструмент.
Например, почти сразу же с его помощью Дж. Спенсер улучшил нижнюю оценку {\it числа Рамсея} (см. определение в [Ga] и задачу \ref{prolov-ram}), которая не поддавалась улучшению в течение сорока лет.
Сейчас диапазон применения леммы становится все шире.
Здесь теория графов и гиперграфов, здесь экстремальные задачи комбинаторики, теория алгоритмов и даже комбинаторная геометрия и теория диофантовых приближений.

За прошедшие десятилетия появились разнообразные усовершенствования Локальной леммы, многие из которых уже лишь отдаленно напоминают первоначальный вариант.
И это еще одно свидетельство исключительной плодотворности идеи Ловаса.

\smallskip
{\it Как устроено изложение в этой заметке.}
Основные идеи демонстрируются по одной и на `олимпиадных' примерах, т.е. на простейших частных случаях,
свободных от технических деталей.
Мы показываем, {\it как можно придумать} лемму Ловаса.
Путь к ее доказательству и применениям намечен в виде задач (всех задач этого и следующего разделов, кроме
задач \ref{prolov-indep} и \ref{prolov-indep105}, которые просто поясняют понятие независимости).
Обучение путем решения задач не только характерно для серьезного изучения математики, но и продолжает древнюю культурную традицию.
\footnote{
Например, послушники дзэнских монастырей обучаются, размышляя над загадками, данными им наставниками.
Впрочем, эти загадки являются скорее парадоксами, а не задачами.
См. подробнее \cite{S}.}

К важнейшим задачам приводятся указания и решения.

Обычно лемму Ловаса излагают на вероятностном языке.
Однако, по нашему мнению, приводимое комбинаторное изложение более доступно и полезно для начинающего.
Важно излагать вероятностные идеи (например, независимости) и развивать вероятностную интуицию, но при этом сохранять строгость изложения.
Разумнее делать это, не определяя понятия вероятностного пространства.
\footnote{Отличие элементарной теории вероятностей от перечислительной комбинаторики скорее в том, что речь идет о {\it долях}  вместо
чисел, и интерес часто представляют {\it оценки}, а не равенства.}
Это как раз подготовит начинающего к введению этого довольно абстрактного понятия,
ср. с [Z, философски-методическое отступление].
Кроме того, вероятностной интуиции начинающего противоречит получение вероятностными методами абсолютно
(а не с некоторой вероятностью) верного результата.
\footnote{Объяснять, как с помощью вероятностных методов можно получить абсолютно верный результат, лучше на более простых примерах.
См., например, задачи \ref{lll-city}, \ref{lll-city3}, \ref{prolov-chain}, \cite[задача 3 на стр. 3]{Go},
\cite[задача 7.3.a]{I}.
Мы хотели бы сделать заметку доступной даже для тех, кто не разбирал таких примеров.}
(Впрочем, для человека, уже владеющего понятием вероятностного пространства, изложение на вероятностном языке не хуже комбинаторного.)

\smallskip
{\it Приведем интересные факты, которые можно доказать при помощи леммы Ловаса (и вряд ли можно доказать  без нее!).}
Видимо, из задач \ref{prob-lovac-2cveta}-\ref{prolov-rainbow} вы сможете решить сейчас только пункты (a).
К пунктам (b) разумно вернуться после изучения следующего раздела.
Более того, задача \ref{prolov-16}.b естественнее по формулировке, но сложнее двух следующих.

\begin{pr}\label{prob-lovac-2cveta}
(a) По каждому из 100 видов работ в фирме имеется ровно 8 специалистов.
Каждому сотруднику нужно дать выходной в субботу или в воскресенье.
Докажите, что это можно сделать так, чтобы и в субботу, и в воскресенье для каждого вида работ присутствовал  специалист по нему.
(Сотрудник может быть специалистом по нескольким видам работ; распределение специалистов по видам работ
известно тому, кто устраивает выходные.)

(b) По каждому из нескольких видов работ в фирме имеется ровно 8 специалистов.
(Теперь видов работ не обязательно 100.)
Каждый вид работ имеет общих специалистов не более, чем с 30 другими видами.
Каждому сотруднику нужно дать выходной в субботу или в воскресенье.
Докажите, что это  можно сделать так, чтобы и в субботу, и в воскресенье для каждого вида работ присутствовал  специалист по нему.
\end{pr}

{\it Замечание.}
Для каждого вида работ $x$ обозначим через $A_x$ множество распределений выходных, при которых
и в субботу, и в воскресенье на работе есть специалист по $x$.
Нужно доказать, что $\cap_x A_x\ne\emptyset$.
В (a) это делается путем подсчета количества элементов.
В (b) этого уже не хватает, нужна идея из следующего раздела.
Там мы покажем, как {\it независимость} (определенную там) можно применять для оценки количества элементов в пересечении множеств.

 Описанную идею можно сформулировать так.
Нужное условие мы представляем в виде пересечения некоторого числа условий.
При этом ясно, что для каждого из них есть конструкция, ему удовлетворяющая.
Иногда отсюда можно вывести, что есть конструкция, удовлетворяющая всем этим условиям одновременно!
Эта идея часто применяется в математике.
(Для читателя, знакомого с соответствующими понятиями, напомним, что в анализе так доказывается существование решения дифференциального уравнения, в топологии --- вложимость $n$-мерного компакта в $\R^{2n+1}$, ср. \cite[\S2]{Sk}.)
Число условий может быть бесконечно, поэтому идея пересечения `равносильна' идее итерационного процесса.
А в настоящей заметке мы покажем, как применять эту идею в комбинаторике.
Несмотря на конечность числа условий, ее применение весьма нетривиально.

\begin{pr}\label{prolov-16}
(a) По кругу стоит 200 студентов из 10 групп, в каждой из которых 20 студентов.
Докажите, что можно в каждой группе выбрать старосту так, чтобы никакие два старосты не стояли рядом.

(b) По кругу стоит 1600 студентов из 100 групп, в каждой из которых 16 студентов.
Докажите, что можно в каждой группе выбрать старосту так, чтобы никакие два старосты не стояли рядом.
 \end{pr}

 \begin{pr}\label{prolov-progr}
(a) Докажите, что можно раскрасить первые 8 натуральных чисел в 2 цвета так, чтобы не было одноцветной арифметической прогрессии длины 3.

(b) Докажите, что можно раскрасить первые 15 миллионов натуральных чисел в 2 цвета так, чтобы не было одноцветной арифметической прогрессии длины 32.
\end{pr}

\begin{pr}\label{prolov-rainbow}
(a) Докажите, что для любого $M\in\R$ можно раскрасить все вещественные числа в 2 цвета так, чтобы для любого $x\in\R$
числа $x$ и $x+M$ были не одноцветны.

(b) Докажите, что для любых 25 чисел $M_1,\dots,M_{25}\in\R$ можно раскрасить все вещественные числа в 3 цвета так, чтобы для любого $x\in\R$
среди чисел $x,x+M_1,\dots,x+M_{25}$ были числа каждого из трех цветов.
\end{pr}

Решения пунктов (b) вышеприведенных задач основаны на идее, аналогичной решению задачи  \ref{prob-lovac-2cveta}.b.

\bigskip
{\bf Независимость и лемма Ловаса. }

Приведем задачи, которые подведут нас к лемме Ловаса \ref{prolov-sym} (почему она интересна, написано в предыдущем разделе).

\begin{pr}\label{lll-city} Каждый житель города либо здоров, либо болен, а также либо богат, либо беден.
Богатство и здоровье {\it независимы}, т.е. доля богатых здоровых среди богатых равна доле здоровых среди всех жителей.
Известно, что есть богатый горожанин и есть здоровый горожанин.
Обязательно ли найдется богатый здоровый горожанин?
\end{pr}

Обозначим через $|X|$ количество элементов в множестве $X$.
Подмножества $A$ и $B$ конечного множества $M$ называются {\it независимыми}, если
$$|A\cap B|\cdot|M| =|A|\cdot|B|.$$
При $B\ne\emptyset$ это равносильно тому, что доля множества $A\cap B$ в $B$ равна доле множества $A$ в $M$.

\begin{pr}\label{prolov-indep}
Зависимы
ли следующие подмножества?
(Мы называем {\it зависимыми} подмножества, не являющиеся независимыми.)

(a) В множестве всех клеток шахматной доски подмножество клеток в первых трех ее строках
и подмножество клеток в последних четырех ее столбцах.

(b) подмножества $\{1,2\}\subset\{1,2,3,4\}$ и $\{1,3\}\subset\{1,2,3,4\}$.

(c) подмножества $\{1,2\}\subset\{1,2,3,4,5,6\}$ и $\{1,3\}\subset\{1,2,3,4,5,6\}$.
\end{pr}

\begin{pr}\label{prolov-indep105}
Зависимы ли следующие подмножества множества целых чисел от 1 до 105?

(a)
Подмножество чисел, делящихся на 5,
и подмножество чисел, делящихся на 7.

(b) Подмножество чисел, делящихся на 15,
и подмножество чисел, делящихся на 21.

(c) Подмножество чисел, делящихся на 15,
и подмножество чисел, делящихся на 5.

(d) Подмножество чисел, делящихся на 10,
и подмножество чисел, делящихся на 7.
\end{pr}

\begin{pr}\label{prolov-indepcol}
(Ср. с замечанием после задачи \ref{prob-lovac-2cveta}.b.)
Зависимы ли следующие подмножества множества всех раскрасок чисел $1,2,\dots,400$ в два цвета?

(a) Подмножество раскрасок, для которых $\{1,2,\dots,8\}$ одноцветно, и  подмножество раскрасок, для которых $\{11,12,\dots,18\}$ одноцветно.

(b) Подмножество раскрасок, для которых $\{1,2,\dots,8\}$ неодноцветно, и  подмножество раскрасок, для которых $\{11,12,\dots,18\}$ неодноцветно
(ср. с задачей \ref{prob-lovac-2cveta}.b).

(c) Подмножество раскрасок, для которых $\{1,2,\dots,8\}$ одноцветно, и  подмножество раскрасок, для которых $\{6,7,\dots,13\}$ одноцветно.
\end{pr}

{\it Если условие задачи является формулировкой утверждения, то в задаче требуется это утверждение доказать.}

\begin{pr}\label{prolov-com}
Подмножества $A$ и $B$ конечного множества независимы тогда и только тогда, когда $A$ и $\overline B$ независимы.
\end{pr}

\begin{pr}\label{lll-city3}
(a) Обязательно ли найдется богатый здоровый умный горожанин, если в городе
доля богатых горожан больше $2/3$, доля  здоровых больше $2/3$, и доля умных больше $2/3$?

(b) Тот же вопрос, если в городе есть богатый горожанин, есть здоровый
горожанин и есть умный горожанин, богатство, здоровье и ум попарно независимы,
и доля богатых здоровых умных среди богатых здоровых такая же, как и доля
умных среди всех жителей.
(Вместе с условием попарной независимости последнее условие называется {\it независимостью в совокупности}.)

(c) Тот же вопрос, если в городе богатых горожан больше половины, здоровых
больше половины, умных больше половины, богатство и ум независимы, здоровье
и ум независимы.
\end{pr}

Задача \ref{lll-city3} показывает, что чем сильнее условие, характеризующее независимость нескольких множеств, тем меньшей доли каждого множества достаточно, чтобы гарантировать непустоту пересечения.
Причем наиболее интересный результат (\ref{lll-city3}.c) получается `посередине' между крайними условиями ---
полного отсутствия независимости (\ref{lll-city3}.a) и независимости в совокупности (\ref{lll-city3}.b).
Так часто бывает: наиболее полезные соображения находятся между `крайними' точками зрения.

\begin{pr}\label{prolov-chain}
(a) Пусть $A_1,A_2,A_3,A_4$ --- подмножества 720-элементного множества, в каждом из которых более 480 элементов.
Если $A_k$ и $A_{k+1}$ независимы для любого $k=1,2,3$, то $A_1\cap A_2\cap A_3\cap A_4\ne\emptyset$.

(b) {\it Гипотеза.} Пусть $n\ge2$ и $A_1,A_2,\dots,A_n$ --- подмножества конечного множества, доля каждого из которых
(т.е. $|A_k|/|M|$) больше $1-\dfrac1{n-1}$.
Если $A_k$ и $A_{k+1}$ независимы для любого $k=1,2,\dots,n-1$, то $A_1\cap A_2\cap \dots\cap A_n\ne\emptyset$.
\end{pr}

Подробнее о независимости см. [KZP].

Для формулировки леммы Ловаса нужно еще более `хитрое' условие независимости на несколько множеств, чем рассмотренные ранее.

Подмножество $A$ конечного множества $M$ называется {\it независимым от набора подмножеств $B_1,\dots,B_k\subset M$}, если $A$ независимо с любым подмножеством, являющимся пересечением нескольких (возможно, одного) множеств из
$B_1,\dots,B_k$.

\begin{pr}\label{prolov-ind} Приведите пример подмножеств $A,B_1,B_2$ конечного множества:

(a) попарно независимых, но для которых $A$ не является независимым от набора $B_1,B_2$;

(b) не являющихся попарно независимыми, но для которых $A$ независимо от набора $B_1,B_2$.
\end{pr}

\begin{pr}\label{prolov-indcol}
Обозначим через $M$ семейство всех раскрасок множества $\{1,2,\dots,400\}$ в два цвета.
Для подмножества $\alpha\subset\{1,2,\dots,400\}$ обозначим через $A_\alpha\subset M$ подмножество тех раскрасок, для которых $\alpha$ одноцветно.
Тогда $A_{\{1,2,\dots,8\}}$ не зависит от набора $\{A_\alpha\ :\ \alpha\subset\{9,10,\dots,400\}\}$.
(Ср. с замечанием после задачи \ref{prob-lovac-2cveta}.b.)
\end{pr}

\begin{pr}\label{prolov-comk} Следующие условия на подмножества $A,B_1,\dots,B_k$ равносильны:

$\bullet$ $A$ независимо от набора $B_1,\dots,B_k$

$\bullet$ $\overline A$ независимо от набора $B_1,\dots,B_k$.

$\bullet$ $\overline A$ независимо от набора $\overline B_1,\dots,\overline B_k$.
\end{pr}

\begin{pr}\label{prolov-sym}
(a) {\bf Локальная лемма Ловаса в симметричной форме.}
{\it Пусть $A_1,\ldots, A_n$ --- подмножества конечного множества.
Если для некоторого $d$ и любого $k$ доля подмножества $A_k$ не меньше $1-\dfrac1{4d}$ и существует набор из не менее,  чем  $n-d$ подмножеств $A_j$, от которого $A_k$ не зависит, то $A_1\cap\dots\cap A_n\ne\emptyset$.}
\footnote{Вот формулировка на вероятностном языке, которая не используется в дальнейшем.
{\it Пусть дано вероятностное пространство и $A_1,\ldots, A_n$ --- события.
Пусть для некоторого $d$ и любого $k$ вероятность события $A_k$ не меньше $1-\frac1{4d}$ и
существует набор из не менее, чем  $n-d$ событий $A_j$, от которого $A_k$ не зависит.
Тогда вероятность события $A_1\cap\ldots\cap A_n$ положительна.}
}

(b) При $d>2$ утверждение пункта (a) верно, если заменить $1-\dfrac1{4d}$ на $1-\dfrac1{e(d+1)}$, где $e$ --- основание натуральных логарифмов.
\end{pr}

Читатель может перед доказательством этой леммы применить ее к решению задачи \ref{prob-lovac-2cveta}.b.
Доказательство леммы нетривиально обобщает идеи решения задач \ref{lll-city3} и \ref{prolov-chain}.
Из этих задач ясно, что нужно оценивать снизу количество элементов в пересечении $s$ из данных множеств, начиная с $s=1$ и заканчивая $s=n$, при помощи индукции по $s$.
Как часто бывает, наиболее трудная часть --- догадаться, какое конкретно утверждение нужно доказывать по индукции
(а также, по каким параметрам вести индукцию).
Вот оно:
$$|A_1\cap\dots\cap A_{k+t}|\ge |A_1\cap\dots\cap A_k|\left(1-\frac1d\right)^t\quad\text{для любых}\quad k,t\ge0.$$

\bigskip
{\bf Указания и решения к задачам \ref{prob-lovac-2cveta}-\ref{prolov-sym} (кроме \ref{prolov-16}.b, \ref{prolov-progr}.b и \ref{prolov-rainbow}.b). }

\smallskip
{\bf \ref{prob-lovac-2cveta}.} (a)
Посчитаем двумя способами количество всех таких пар $(a,x)$, что $a$ --- распределение выходных и $x$ --- вид работы, по которому в один из дней не будет специалиста при распределении $a$.
Обозначим через $n$ общее число специалистов.
Для каждого вида работ имеется $2\cdot2^{n-8}=2^{n-7}$ распределений выходных, при которых все специалисты по этому
виду работ отдыхают в один и тот же день.
Так как видов работ 100, то количество пар не больше $100\cdot 2^{n-7}<2^n$.
Общее число распределений выходных равно $2^n$.
Значит, найдётся распределение выходных, при котором для каждого вида работ не все специалисты по нему отдыхают
в один и тот же день.

(b) Обозначим через $A$ множество распределений выходных.
Для каждого вида работ $x$ обозначим через $\widehat x$ множество специалистов по нему, а через $A_x$ --- множество распределений выходных, при которых и в субботу, и в воскресенье на работе есть специалист по $x$.
Тогда $|A_x|/|A|=2^{-7}$.
Подмножество $A_x$ не зависит от набора $\{A_y\ |\ \widehat y\cap \widehat x=\varnothing\}$.
Так как каждый вид работ имеет общих специалистов не более, чем с 30 другими видами, то вне этого набора не более 30 подмножеств.
Применим локальную лемму Ловаса в симметричной форме (задачу \ref{prolov-sym}.a)
к дополнениям множеств $A_x$ и $d=2^5$.
Это возможно ввиду утверждения задачи \ref{prolov-com} и неравенства $30<2^5$.
Получим $\cap_x\overline{A_x}\ne\varnothing$.

\smallskip
{\bf \ref{prolov-16}.} (a) Выберем произвольного студента произвольной группы и назначим его старостой.
Далее действуем так: на каждом шаге выбираем из группы, в которой еще не выбран староста, студента, не являющегося соседом никакого выбранного старосты.
Так как выбранных старост не больше $9$, то соседей у выбранных старост не больше $18$.
Следовательно, на каждом шаге мы можем найти нужного студента.
В конце получим $10$ человек из разных групп, никакие два из которых не являются соседями.


\smallskip
{\bf \ref{prolov-rainbow}.} (a) Покрасим каждое число $x\in\R$ в четность числа $[x/M]$.

\smallskip
{\bf \ref{prolov-indep}.} Ответы: (a,b) независимы, (c) зависимы.

\smallskip
{\bf \ref{prolov-indep105}.} Ответы: (a) независимы, (b,c,d) зависимы.

\smallskip
{\bf \ref{prolov-indepcol}.} Ответы: (a,b) независимы, (c) зависимы.

\smallskip
{\bf \ref{lll-city3}.} (c)
Забудьте про глупых людей!

Приведем более сложное решение.
Зато оно подводит к лемме Ловаса.
Обозначим через $\text{У,Б,З}$ множества умных, богатых и здоровых горожан.
Будем пропускать знаки пересечения и числа элементов.
Тогда $\text{УБ}>\text{У}/2<\text{УЗ}$.
Значит,
$$\text{УБЗ}=\text{УБ}-\text{УБ}\overline{\text{З}}>\frac{\text{У}}2-\text{У}\overline{\text{З}}>(\frac12-1+\frac12)\text{У}=0.$$

\smallskip
{\bf \ref{prolov-chain}.}
(a) Будем пропускать знаки пересечения и числа элементов.
Не уменьшая общности, $A_2\ge A_3$.
Тогда аналогично решению задачи \ref{lll-city3}.c $A_1A_2A_3>A_2/3$.
Поэтому
$$A_1A_2A_3A_4=A_1A_2A_3-A_1A_2A_3\overline{A_4}>(A_2/3)-A_3\overline{A_4}>(A_2-A_3)/3\ge0.$$

\comment


(b) Утверждение вытекает из следующего результата для $p=1- \frac{1}{k-1}$.

{\it Пусть $A_1,A_2,\dots,A_n$ --- подмножества конечного множества, доля каждого из которых больше $p$.
Если $A_k$ и $A_{k+1}$ независимы для любого $k=1,2,\dots,n-1$, то $|A_1\cap A_2\cap \dots\cap A_n|> ((n-1)p^2-(n-2)p)|M|$.}

Докажем этот результат по индукции.
База индукции $n=1$ следует из $\frac{|A_1|}{|M|}>p$.

Для доказательства шага индукции  будем пропускать знаки пересечения и числа элементов.
Тогда
$$A_1\dots A_n=A_1\dots A_{n-1}-A_1\dots A_{n-1}\overline A_n \ge A_1\dots A_{n-1}-A_{n-1}\overline A_n>$$ $$>((n-2)p^2-(n-3)p)M-(1-p)pM=((n-1)p^2-(n-2)p)M.$$



\endcomment

{\bf \ref{prolov-sym}.}
(a) Можно считать, что $d>1$, иначе утверждение очевидно.
Достаточно доказать для любых $k,t\ge0$ неравенство
$$I(k+t,t):\qquad |A_1\cap\ldots\cap A_{k+t}|\ge |A_1\cap\ldots\cap A_k|\left(1-\frac1d\right)^t.$$
Утверждение леммы получается из $I(n,n)$.
Докажем неравенство $I(k+t,t)$ индукцией по паре
\footnote{Подумайте, почему не проходит индукция по паре $(k,t)$ или по $s+t=k+2t$, а также зачем рассматривать $k>0$, если лемма Ловаса вытекает из случая $k=0$.}
$(s,t)$, где $s=k+t$.
База $s=k+t=t=0$ очевидна.

Докажем шаг, т.е. докажем неравенство $I(k+t,t)$, предполагая выполненным неравенство $I(k'+t',t')$ для всех пар $(k'+t',t')$, лексикографически меньших, чем пара $(k+t,t)$.
Тогда $k+t>0$.
Случай $t=0$ очевиден.
Пусть теперь $t>0$.
Oбозначим $A_{1,2,\dots, m}:=A_1\cap A_2\cap\ldots\cap A_m$.

{\it Случай $t>1$} сводится к последовательному применению
неравенств $I(k+t,1)$ и $I(k+t-1,t-1)$:
$$|A_{1,2,\dots,k+t}|\ge |A_{1,2,\dots,k+t-1}|\left(1-\frac1d\right)\ge |A_{1,2,\dots,k}|\left(1-\frac1d\right)^{1+t-1}=
|A_{1,2,\dots,k}|\left(1-\frac1d\right)^t.$$

{\it Пусть $t=1$. }
По условию существует $X\subset \{1,2,\dots,n\}$ такое, что $|X|\ge n-d$ и $A_{k+1}$ не зависит от набора
$\{A_j\ :\ j\in X\}$.
Тогда
$$|\{1,2,\dots,k\}-X|\le n-|X|\le n-(n-d)=d.$$
Поэтому можно считать, что $A_{k+1}$ не зависит от набора $A_1,A_2,\dots,A_{k-d}$ (при $k\le d$ этот набор пуст).
Тогда
$$|\overline{A_{k+1}}\cap A_{1,2,\dots,k}|\overset{(1)}\le
|\overline{A_{k+1}}\cap A_{1,2,\dots,k-d}|\overset{(2)}\le
\frac{|A_{1,2,\dots,k-d}|}{4d}\overset{(3)}\le
\frac{|A_{1,2,\dots,k}|}{4d(1-\frac1d)^d}\overset{(4)}\le
\frac{|A_{1,2,\dots,k}|}d.$$
Здесь

$\bullet$ при $k\le d$ по определению считаем, что $A_{1,2,\dots,k-d}$ есть то данное в условии множество,
подмножествами которого являются $A_1,\dots,A_n$;

$\bullet$ первое неравенство справедливо ввиду $A_{1,2,\dots,k}\subset A_{1,2,\dots,k-d}$;

$\bullet$ второе --- ввиду того, что доля подмножества $A_{k+1}$ не меньше $1-\dfrac1{4d}$ и, при $k>d$,
ввиду независимости $\overline{A_{k+1}}$ от набора $A_1,A_2,\dots,A_{k-d}$ (см. задачу \ref{prolov-com});

$\bullet$ третье есть $I(k,d)$ и верно по предположению индукции;

$\bullet$ четвертое --- ввиду $\left(1-\dfrac1d\right)^d\ge\left(1-\dfrac12\right)^2=\dfrac14$.

Значит,
$$|A_{1,2,\dots,k+1}|=|A_{1,2,\dots,k}|-|\overline{A_{k+1}}\cap A_{1,2,\dots,k}|\ge |A_{1,2,\dots,k}|\left(1-\dfrac1d\right).$$

\comment

\smallskip
{\it Другое изложение этого доказательства.}
Можно считать, что $d>1$, иначе утверждение очевидно.
Достаточно доказать индукцией по $k$, что
$$|A_1\cap\dots\cap A_{k+1}|\ge |A_1\cap\dots\cap A_k|\left(1-\frac1d\right)\quad\text{для любого}\quad k,\ 0\le k\le n-1.$$
Перемножив все эти неравенства, получим утверждение леммы (аккуратно это доказываетя по индукции, см. случай $t>1$ в предыдущем изложении доказательства).

База $k=0$ очевидна. Докажем шаг. Положим $A_{1,2,\dots, m}:=A_1\cap A_2\cap\ldots\cap A_m$ для любого $m$.

По условию существует $X\subset \{1,2,\dots,n\}$ такое, что $|X|\ge n-d$ и $A_{k+1}$ не зависит от набора
$\{A_j\ :\ j\in X\}$.
Тогда
$$|\{1,2,\dots,k\}-X|\le n-|X|\le n-(n-d)=d.$$
Поэтому можно считать, что $A_{k+1}$ не зависит от набора $A_1,A_2,\dots,A_{k-d}$ (при $k\le d$ этот набор пуст).
Тогда
$$|\overline{A_{k+1}}\cap A_{1,2,\dots,k}|\overset{(1)}\le
|\overline{A_{k+1}}\cap A_{1,2,\dots,k-d}|\overset{(2)}\le
\frac{|A_{1,2,\dots,k-d}|}{4d}\overset{(3)}\le
\frac{|A_{1,2,\dots,k}|}{4d(1-\frac1d)^d}\overset{(4)}\le
\frac{|A_{1,2,\dots,k}|}d.$$
Здесь

$\bullet$ при $k\le d$ в качестве $A_{1,2,\dots,k-d}$ следует взять
множество из условия леммы, подмножествами которого являются $A_1,\dots,A_n$;

$\bullet$ первое неравенство справедливо ввиду $A_{1,2,\dots,k}\subset A_{1,2,\dots,k-d}$;

$\bullet$ второе --- ввиду того, что доля подмножества $A_{k+1}$ не меньше $1-\dfrac1{4d}$ и, при $k>d$,
ввиду независимости $\overline{A_{k+1}}$ от набора $A_1,A_2,\dots,A_{k-d}$ (см. задачу
9);

$\bullet$ третье ---  по предположению индукции: перемножим доказанн!!!ые неравенства для $k-d,\dots,k-1$
(аккуратно это доказываетя по индукции, см. случай $t>1$ в предыдущем изложении доказательства);

$\bullet$ четвертое --- ввиду $\left(1-\dfrac1d\right)^d\ge\left(1-\dfrac12\right)^2=\dfrac14$.

Значит,
$$|A_{1,2,\dots,k+1}|=|A_{1,2,\dots,k}|-|\overline{A_{k+1}}\cap A_{1,2,\dots,k}|\ge |A_{1,2,\dots,k}|\left(1-\dfrac1d\right).$$

\endcomment

\bigskip
{\bf Указания и решения  к задачам \ref{prolov-16}.b, \ref{prolov-progr}.b и \ref{prolov-rainbow}.b.}

\smallskip
{\bf \ref{prolov-16}.}
(b) Обозначим через $A$ семейство подмножеств из 100 студентов, в которых по одному студенту из каждой группы.
Для любого студента $x$ обозначим через $A_x$ семейство подмножеств из $A$, содержащих и студента $x$,
и следующего за ним по часовой стрелке студента $x_+$.

Если $x$ и $x_+$ одногруппники, то $A_x$ пусто.
Иначе $|A_x|/|A| = 1/256$.
Итак, всегда $|A_x|/|A| \leq 1/256$.

Множество $A_x$ не зависит от набора $\alpha_x$ всех тех $A_y$, для которых ни один из студентов $y,y_+$ не является одногруппником ни со студентом $x$, ни со студентом $x_+$.
Семейство $A_z$ не входит в набор $\alpha_x$ тогда и только тогда, когда один из студентов $z,z_+$ является одногруппником с одним из студентов $x,x_+$.
Тех $z$, которые являются одногруппниками для $x$, ровно 16.
Аналогичное верно с заменой пары $(x,z)$ на любую из пар $(x_+,z)$, $(x,z_+)$, $(x_+,z_+)$.
Поэтому количество студентов $z$, для которых $A_z\not\in\alpha_x$, не больше  $16\cdot4=64$.

Применим локальную лемму Ловаса в симметричной форме (задачу \ref{prolov-sym}.a)
к дополнениям множеств $A_x$ и  $d=64$.
Это возможно ввиду утверждения задачи \ref{prolov-com} и равенства $4\cdot 64=256$.
Получим $\cap_x\overline{A_x}\ne\emptyset$.

\smallskip
{\bf \ref{prolov-progr}.}
(b) 
Положим
$n:=15\cdot10^6$.
Обозначим через $A$ семейство раскрасок множества $\{1,2,\dots,n\}$ в 2 цвета.
Для любой 32-элементной арифметической прогрессии $\alpha\subset\{1,2,\dots,n\}$ обозначим через $A_\alpha$ семейство раскрасок множества $\{1,2,\dots,n\}$ в 2 цвета, для которых $\alpha$ одноцветна.
Тогда $|A_\alpha|/|A|=2^{-31}$.
Подмножество $A_\alpha$ не зависит от набора $\{A_\beta \ |\ \beta\cap\alpha=\emptyset\}.$
Каждую 32-элементную арифметическую прогрессию в $\{1,2,\dots,n\}$ пересекает не более $32^2[n/31]<2^{10}2^{19}=2^{29}$ других таких прогрессий (докажите!).
Применим локальную лемму Ловаса в симметричной форме (задачу \ref{prolov-sym}.a)
к дополнениям множеств $A_\alpha$ и  $d=2^{29}$.
Это возможно ввиду утверждения задачи \ref{prolov-com}.
Получим $\cap_\alpha\overline{A_\alpha}\ne\emptyset$.

\smallskip
{\bf \ref{prolov-rainbow}.}
(b) Докажем следующее более слабое утверждение.

{\it Для любых 25 чисел $M_1,\dots,M_{25}\in\R$ и конечного множества $X\subset \R$ можно раскрасить все вещественные числа в 3 цвета
так, чтобы для любого $x\in X$ среди чисел $x,x+M_1,\dots,x+M_{25}$ были числа каждого цвета.}

При помощи {\it соображений компактности} можно вывести из этого утверждения его аналог для бесконечного $X$, в частности, для $X=\R$ [AS, \S5.2].

Обозначим через $A$ семейство раскрасок множества $X\cup(M_1+X)\cup\dots\cup(M_{25}+X)$ в 3 цвета.
Для любого $x\in X$ обозначим через $A_x\subset A$ подсемейство раскрасок, для которых среди цветов чисел $x,x+M_1,\dots,x+M_{25}$ не все цвета присутствуют.
Тогда $|A_x|/|A|\le3(2/3)^{26}$.
Каждое множество $A_x$ `зависимо не более, чем с $25\cdot26=650$ другими' (т.е. независимо от набора всех множеств, кроме некоторых 650).
Применим локальную лемму Ловаса в симметричной форме (задачу \ref{prolov-sym}.a)
к дополнениям множеств $A_x$ и  $d=650$.
Это возможно ввиду утверждения задачи \ref{prolov-com} и неравенства $(3/2)^{26}>2^{13}>8000>7800=3\cdot2600=3\cdot4\cdot650$.
Получим $\cap_x\overline{A_x}\ne\emptyset$.

 \bigskip
{\bf Задачи для самостоятельного решения. }

Указания и решения к этим задачам можно найти в [AS, I].

\begin{pr}\label{prolov-hyp}
Дано число

(a) $k\ge 10$; \quad (b) $k=9$

и семейство $k$-элементных подмножеств конечного множества $M$.
Если каждый элемент множества $M$ содержится ровно в $k$ подмножествах семейства, то существует раскраска
множества $M$ в два цвета, для которой каждое подмножество семейства содержит элементы обоих цветов.
(Т.е. хроматическое число любого $k$-однородного $k$-регулярного гиперграфа равно двум при $k\ge9$.
Ср. с задачей \ref{prob-lovac-2cveta}.b.)
\end{pr}

\begin{pr}\label{prolov-subsets}
В конечном множестве выбрано несколько подмножеств.
В каждом из них не менее 3 элементов.
Каждое из них пересекается не более чем с $a_{i}$ выбранными $i$-элементными подмножествами.
Если $\sum_i a_i2^{-i} \leq 1/8$,
то можно покрасить элементы данного множества в два цвета так, чтобы каждое выбранное подмножество содержало элементы обоих цветов.
\end{pr}

\begin{pr}\label{prob-lovac-cikl} (a) Для любого разбиения множества вершин цикла длины $11n$ на $n$ множеств
по 11 вершин можно выбрать по вершине из каждого множества так, что между выбранными $n$ вершинами нет ребер.


(b) В графе степень каждой вершины не превосходит $\Delta$.
Все вершины раскрашены в $r$ цветов.
Вершин каждого цвета не менее $2e\Delta + 1$.
Тогда можно выбрать $r$ вершин разных цветов, никакие две из которых не соединены ребром.
\end{pr}

\begin{pr}\label{prob-lovac-wdw}
(a)  Каждую $k$-элементную арифметическую прогрессию в множестве $\{1,2,\dots,n\}$ пересекает не более
$k^2[\dfrac n{k-1}]$ других таких прогрессий.

(b) Для любого натурального $k$ существует раскраска первых
$[2^{k-3}(k-1)/k^2]$ натуральных чисел в 2 цвета,
для которой нет одноцветной $k$-элементной арифметической прогрессии.

(c) Каждую $k$-элементную арифметическую прогрессию в множестве $\{1,2,\dots,n\}$ пересекает не более
$nk$ других таких прогрессий.

(d) Для любого натурального $k$ существует раскраска первых $[2^{k-3}/k]$ натуральных чисел в 2 цвета,
для которой нет одноцветной $k$-элементной арифметической прогрессии.
\end{pr}

\begin{pr}\label{prob-lovac-raduga}
(a) Если $X\subset \R$ --- конечное множество  и $m,r$ --- натуральные числа, для которых
$4rm(m-1)\left(1-\dfrac 1 r\right)^m<1,$
то для любого $m$-элементного подмножества $M\subset\R$ существует раскраска множества $\R$ в $r$ цветов такая, что
для любого $x\in X$ множество $x+M:=\{x+a\ :\ a\in M\}$
содержит точки каждого из $r$ цветов.

(b) То же для $X=\Z$.

(c) То же для $X=\R$.
\end{pr}

\begin{pr}\label{prolov-ram}
Определение {\it числа Рамсея}  $R(n,n)$ см., например, в [Ga].

(a) Если $\displaystyle{n\choose 2}{k\choose n-2}+1<2^{{n\choose 2}-1}/e$, то $R(n,n)>k$,

(b) $R(n,n)>\sqrt2e^{-1}n2^{n/2}(1+o(1))$.
\end{pr}

\begin{pr}\label{prob-lovas-10d} Имеется несколько цветов.
Каждой вершине некоторого графа сопоставлен список из не менее, чем $10d$ этих цветов, где $d>1$.
Для любых вершины $v$ и цвета из ее списка имеется не более $d$ соседей вершины $v$, в списке которых есть этот цвет.  Докажите, что можно так раскрасить каждую вершину графа в некоторый цвет из ее списка, чтобы концы любого ребра были разноцветны.
\end{pr}

\begin{pr}\label{prolov-orgr}
В ориентированном графе
в каждую вершину входит не больше $\Delta$ ребер и из каждой вершины выходит не меньше $\delta$ ребер.
Тогда для любого натурального $k \le \dfrac1{1 - (4\delta \Delta)^{-1/\delta}}$ в графе найдется ориентированный цикл длины, кратной $k$.
\end{pr}

\begin{pr}\label{prolov-latin}
Клетки доски $n \times n$ раскрашены в несколько цветов.
Клеток каждого цвета не больше, чем $(n-1)/16$.
Тогда можно поставить на доску $n$ попарно не бьющих друг друга ладей, чтобы они стояли на клетках разных цветов.
\end{pr}

\begin{pr}\label{prob-lovac-SAT}
{\it КНФ-формула} или {\it формула в конъюнктивной нормальной форме} ---  конъюнкция набора дизъюнкций нескольких из переменных $x_1,\dots,x_n$ и их отрицаний.
 Если в каждом `сомножителе' КНФ-формулы ровно $k$ `слагаемых' и у каждого `сомножителя' есть общие переменные не более, чем с $2^{k-2}$ другими, то булева функция, определяемая формулой, не является тождественным нулем.
\end{pr}

{\it Замечание.}
Одной из центральных в информатике является проблема $k$-выполнимости ($k$-SAT problem):
{\it существует ли алгоритм, который по КНФ-формуле, в каждой дизъюнкции которой ровно $k$ переменных, выясняет, является ли тождественным нулем
соответствующая булева функция.}
 При $k=2$ есть полиномиальный алгоритм её решения.
При б\'ольших $k$ это наиболее стандартная {\it NP-полная проблема}.
Поэтому полиномиальный алгоритм дал бы и равенство {\it классов P и NP}.

\begin{pr}\label{prolov-gen}
(a) {\bf Локальная лемма Ловаса.}
{\it Пусть $A_1,\ldots, A_n$ --- подмножества конечного множества, $J_1,\ldots,J_n\subset\{1,\ldots,n\}$ и $\gamma_1,\dots,\gamma_n\in(0,1)$.
Пусть для любого $k$

$\bullet$ доля подмножества $A_k$ не меньше $1-(1-\gamma_k)\prod_{j\not\in J_k}\gamma_j$;

$\bullet$ множество $A_k$ не зависит от набора $\{A_j\ :\ j\in J_k\}$.

Тогда доля пересечения $\bigcap_{k=1}^nA_k$ не меньше $\prod_{k=1}^n\gamma_k>0$.}
\footnote{Вот формулировка на вероятностном языке.
{\it Пусть дано вероятностное пространство, $A_1,\ldots, A_n$ --- события, $J_1,\ldots,J_n\subset\{1,\ldots,n\}$ и
$\gamma_1,\dots,\gamma_n\in(0,1)$.
Пусть для любого $k$ вероятность события $A_k$ не меньше $1-(1-\gamma_k)\prod_{j\not\in J_k}\gamma_j$ и
событие $A_k$ не зависит от набора $\{A_j, j\in J_k\}$.
Тогда вероятность события $A_1\cap\ldots\cap A_n$ не меньше $\prod_{j=1}^m\gamma_j$.}
}

(b) Существует такое $c>0$, что $R(3,n)>cn\sqrt n$ для любого $n$.
 \end{pr}



\end{document}